\documentclass[12pt]{article}

\usepackage{amsmath}
\usepackage{amsthm}
\usepackage{amssymb}
\usepackage{breqn}
\usepackage{graphicx}
\usepackage{float}

\theoremstyle{plain}
\newtheorem{theorem}{Theorem}[section]
\newtheorem{definition}[theorem]{Definition}

\newtheorem{property}[theorem]{Property}
\newtheorem{lemma}[theorem]{Lemma}
\newtheorem{corollary}[theorem]{Corollary}
\newtheorem{example}[theorem]{Example}

\newcommand{\generateur}[1]{#1_0}
\newcommand{\prefix}[1]{#1_1}
\newcommand{\suffix}[1]{#1_2}

\DeclareMathOperator{\lcp}{lcp}
\DeclareMathOperator{\lcs}{lcs}

\usepackage{color}

\title{Combinatorics of the Interrupted Period.}
\author{A. Thierry}
\begin{document}
\maketitle

\begin{abstract}
This article is about discrete periodicities and their combinatorial structures. It presents and describes the unique structure caused by the alteration of a pattern in a repetition. Those alterations of a pattern arise in the context of double squares and were discovered while working on bounding the number of distinct squares in a string. Nevertheless, they can arise in other phenomena and are worth being presented on their own.
\end{abstract}


\bigskip

If $\generateur{x}$ is a primitive word, and $\prefix{x}$ a prefix of $\generateur{x}$, the sequence $\generateur{x}^n\prefix{x}\generateur{x}^m$ has a singularity: it has a periodic part of period $\generateur{x}$, an interruption, and a resumption of the pattern $\generateur{x}$. That interruption creates a different pattern, one that does not appear in $\generateur{x}^n$. The goal of this article is to unveil that pattern.

\section{Preliminaries}

In this section, we introduce the notations and present a simple property and two of its corollaries. These observations are straightforward but their proofs introduce the
technique used to prove Theorem \ref{ath} and provide insights.
   
We first fix some notations.
 An \emph{alphabet} $A$ is a finite set. 
 We call \emph{letters} the elements of $A$. If $\lvert A \rvert = 2$, the words are referred to as binary and are used in computers. Another well known example for $\lvert A \rvert = 4$ is DNA. \\
 A vector of $A^n$ is a \emph{word} $w$ of length $\lvert w \rvert = n$, which can also be presented under the form of an array $w[1,...,n]$.
Two words are \emph{homographic} if they are equal to each other.
If $x = \prefix{x}\suffix{x}x_3$ for non-empty words $\prefix{x}, \suffix{x}$ and $x_3$, then $\prefix{x}$ is a \emph{prefix} of $x$, $\suffix{x}$ is a \emph{factor} of $x$, and $x_3$ is a \emph{suffix} of $x$ (if both the prefix and the suffix are non empty, we refer to them as proper). We define \emph{multiplication} as concatenation. In english, $breakfast = break . fast$.
In a traditional fashion, we define the \emph{$n^{th}$ power} of a word $w$ as $n$ time the multiplication of $w$ with itself.
 A word $x$ is \emph{primitive} if $x$ cannot be expressed as a non-trivial power of another word $x'$.\\
 A word $\tilde{x}$ is a \emph{conjugate} of $x$ if $x=\prefix{x}\suffix{x}$ and $\tilde{x}=\suffix{x}\prefix{x}$ for non-empty words $\prefix{x}$ and $\suffix{x}$. The set of conjugates of $x$ together with $x$ form the conjugacy class of $x$ which is denoted $Cl(x)$. \\
 A factor $x, \lvert x \rvert =n$ of $w$ has \emph{period} $p$ if $x[i]=x[i+\lvert p \rvert], \forall i \in [1,...n-\lvert p\rvert]$.\\
The \emph{number of occurrences} of a letter $c$ in a word $w$ is denoted $n_c(w)$,  the \emph{longest common prefix} of $x$ and $y$ as $lcp(x,y)$ , while $lcs(x,y)$ denotes the \emph{longest common suffix} of $x$ and $y$ (note that $\lcs (x,y)$ and $\lcp (x,y)$ are words).\\\\


The properties presented next rely on a simple counting argument. If the proofs are not interesting in themselves, they still allow for meaningful results.

\begin{property}\label{fr}
A word $w$ and all of its conjugates have the same number of occurrences for all of their letters, i.e. $\forall\tilde{w} \in Cl(w), \forall a \in A,\ n_a(w) = n_a(\tilde{w})$.
\end{property}

\begin{proof}
Note that $\forall \tilde{w} \in Cl(w), \exists w_1, w_2$, such that $w = w_1w_2, \tilde{w}=w_2w_1$. Then, $\forall a \in A, n_a(w) = n_a(w_1) + n_a(w_2) = n_a(\tilde{w})$. \qed
\end{proof}


The negation of Property \ref{fr} gives the following corollary:

\begin{corollary}\label{number}
If two words do not have the same number of occurrence for the same letter, they are not conjugates.
\end{corollary}

Another important corollary of Property \ref{fr} is the following:

\begin{corollary}
Let $x$ be a word, $\lvert x \rvert \geq n+1$. If $u=x[1...n]$ and $v=x[2...n+1]$ are conjugates of each other, then $x[1] = x[n+1]$, i.e. $v$ is a cyclic shift of $u$.
\end{corollary}

\begin{proof}
Note that $u$ and $v$ have the factor $x[2...n]$ in common. Since $u$ and $v$ are conjugates, they have the same number of occurrences for all of their letters (Proposition \ref{fr}). It follows that $n_{x[1]}(u) = n_{x[1]}(x[1...n]) = n_{x[1]}(x[2...n]) + 1  = n_{x[1]}(v) = n_{x[1]}(x[2...n]) +n_{x[1]}(x[n+1])$, hence $n_{x[1]}(x[n+1]) = 1$, i.e. $x[1] = x[n+1]$. \qed
\end{proof}
\section{Theorem}

  Discrete periods were described by N.J. Fine and H.S. Wilf in 1965 in the article ``Uniqueness theorem for periodic functions'' \cite{FW65}. A corollary of that theorem, the synchronization principle, was proved by W. Smyth in \cite{S05} and L. Ilie in \cite{I05}:
  
  \begin{theorem}\label{fw}
  If $w$ is primitive, then, for all conjugates $\tilde{w}$ of $w, w \neq \tilde{w}$.
  \end{theorem}
  
Which is about the synchronization of patterns. The next theorem is about the impossible synchronization when a pattern is interrupted.
  

First, we need to formalize what we call an interruption of the pattern. Let $\generateur{x}$ be a primitive word and $\prefix{x}$ be a proper prefix of $\generateur{x}$, i.e. $\prefix{x}\ne \generateur{x}$. Write $\generateur{x}=\prefix{x}\suffix{x}$ for some suffix $\suffix{x}$ of $\generateur{x}$. \\\\
Let $W=\generateur{x}^{e_1}\prefix{x}\generateur{x}^{e_2}$ with $e_1\geq1, e_2 \geq 1, e_1+e_2\geq 3$.\\\\
We see that $W$ has a repetition of a pattern $\generateur{x}$ as a prefix: $\generateur{x}^{e_1}\prefix{x}$, and then the repetition is interrupted at position $\lvert \generateur{x}^{e_1}\prefix{x} \rvert$, before starting again in the suffix $\generateur{x}^{e_2}$.
We need one more definition (albeit that definition is not necessary, it is presented here for better understanding) before introducing the two factors that we claim have restricted occurrences in $W$. 

\begin{definition}
Let $W=\generateur{x}^{e_1}\prefix{x}\generateur{x}^{e_2}$ with $e_1\geq1, e_2 \geq 1, e_1+e_2\geq 3$ for a primitive word $\generateur{x}=\prefix{x}\suffix{x}$. Let $\tilde{p}$ be the prefix of length $\lvert \lcp (\prefix{x}\suffix{x}, \suffix{x}\prefix{x})\rvert+1$ of $\prefix{x}\suffix{x}$ and $\tilde{s}$ the suffix of length $\lvert \lcs (\prefix{x}\suffix{x}, \suffix{x}\prefix{x})\rvert+1$ of $\suffix{x}\prefix{x}$. The factor $\tilde{s}\tilde{p}$ starting at position $\lvert \generateur{x}^{e_1} \rvert + \lvert \prefix{x} \rvert - \lvert \lcs (\prefix{x}\suffix{x}, \suffix{x}\prefix{x})\rvert -1$ is the \emph{core of the interrupt} of $W$.\\
\end{definition}
If $W$ and its interrupt are clear from the context, we will just speak of the core (of the interrupt).

\bigskip

\begin{example}
Consider $\generateur{x} = aaabaaaaaabaaaa$ and $\prefix{x} = aaabaaaaaabaaa$, then $\generateur{x}\prefix{x}\generateur{x}^2$ has $\generateur{x}\prefix{x}\generateur{x} = aaabaaaaaabaaaa\mathbf{aaabaaaaaabaaa}aaabaaaaaabaaaa$ as a prefix and $\suffix{x} =  a$. It follows that $\lcp (\prefix{x}\suffix{x},\suffix{x}\prefix{x}) = aaa$, and $ \tilde{p} = aaab, \lcs (\prefix{x}\suffix{x},\suffix{x}\prefix{x}) = aaa$, and $\tilde{s}=baaa$. The core of the interrupt, $\tilde{s}\tilde{p}$, is  underlined in: \[ x\prefix{x}x = aaabaaaaaabaaaaaaabaaaaaa\underbrace{baaaaaab}_{\tilde{s}\tilde{p}}aaaaaabaaaa. \]
\end{example}

The factors that were previously known to have restricted occurrences in $W$, to the best of the author's knowledge, were the inversion factors defined by A. Deza, F. Franek and A. Thierry in \cite{DFT15}:

\begin{definition}
Let $W=\generateur{x}^{e_1}\prefix{x}\generateur{x}^{e_2}$ with $\generateur{x}=\prefix{x}\suffix{x}$ a primitive word and $e_1\geq1, e_2 \geq 1, e_1+e_2\geq 3$. An \emph{inversion factor} of $W$ is a factor that starts at position $i$ and for which:
\begin{itemize}
\item $W[i+j]=W[i+j+\lvert \generateur{x} \rvert + \lvert \prefix{x}\rvert]$ for $0\leq j < \lvert \prefix{x}\rvert$, and
\item $W[i+j]=W[i+j+\lvert \prefix{x}\rvert ]$ for $ \lvert \prefix{x}\rvert \leq j \leq  \lvert x\rvert +  \lvert \prefix{x}\rvert$.
\end{itemize}
\end{definition}

Those inversion factors, which have the structure of $\suffix{x}\prefix{x}\prefix{x}\suffix{x} = \tilde{\generateur{x}}\generateur{x}$, and which length are twice the length of $\generateur{x}$, were used as two notches that forces a certain synchronization of certain squares in the problem of the maximal number of squares in a word, and allowed to offer a new bound to that problem. The main anticipated application of the next result is an improvement of that bound, though the technique has already proved useful in the improvement of M. Crochemore and W. Rytter's three squares lemma, \cite{CR95}, by H. Bay, A. Deza and F. Franek, \cite{BDF15}, and in the proof of the New Periodicity Lemma by H. Bay, F. Franek and W. Smyth \cite{BFS15}.

\bigskip

Now, let $w_1$ be the factor of length $\lvert \generateur{x} \rvert$ of $W$ that has the core of the interrupt of $W$ as a suffix,  and let $w_2$ be the factor of length $\lvert x \rvert$ that has the core of the interrupt of $W$ as a prefix. We will show that both $w_1$ and $w_2$ have restricted occurrences in $W$.\\

\begin{theorem}\label{ath}
Let $\generateur{x}$ be a primitive word, $\prefix{x}$ a proper prefix of $\generateur{x}$ and $W=\generateur{x}^{e_1}\prefix{x}\generateur{x}^{e_2}$ with $e_1\geq1, e_2 \geq 1, e_1+e_2\geq 3$. Let $w_1$ be the factor of length $\lvert \generateur{x} \rvert$ of $W$ ending with the core of the interrupt of $W$, and let $w_2$ be the factor of length $\lvert \generateur{x} \rvert$ starting with the core of the interrupt of $W$. The words $w_1$ and $w_2$ are not in the conjugacy class of $\generateur{x}$.\\
\end{theorem}

\begin{proof}
Define $p=\lcp(\prefix{x}\suffix{x}, \suffix{x}\prefix{x})$ and $s=\lcs(\prefix{x}\suffix{x}, \suffix{x}\prefix{x})$ (note that $p$ and $s$ can be empty). \\
Deza, Franek, and Thierry  showed that $\lvert \lcs(\prefix{x}\suffix{x},\suffix{x}\prefix{x})\rvert+\lvert \lcp(\prefix{x}\suffix{x},\suffix{x}\prefix{x})\rvert\leq\lvert \prefix{x}\suffix{x}\rvert-2$ when $\prefix{x}\suffix{x}$ is primitive (see \cite{DFT15}). Note that in the case $\lvert \lcs(\prefix{x}\suffix{x},\suffix{x}\prefix{x})\rvert+\lvert \lcp(\prefix{x}\suffix{x},\suffix{x}\prefix{x})\rvert = \lvert x\rvert-2$, $w_1$ $w_2$ are the same factor.   \\
Write $\generateur{x}=pr_prr_ss$ and $\tilde{\generateur{x}}=pr'_pr'r'_ss$ for the letters $r_p, r'_p, r_s, r'_s$, $r_p \neq r'_p, r_s \neq r'_s$ (by maximality of the longest common prefix and suffix) and the possibly empty and possibly homographic words $r$ and $r'$.\\
We have, by construction, $w_1=r'r'_sspr_p$ and $w_2=r'_sspr_pr$.\\
Note that $n_{r_p}(w_1) = n_{r_p}(\tilde{\generateur{x}}) + 1$ and that $n_{r'_p}(\tilde{\generateur{x}}) = n_{r'_p}(w_1) + 1$ and, by Corollary \ref{number}, $w_1$ is not a conjugate of $\tilde{\generateur{x}}$, nor of $\generateur{x}$. And because $\lvert w_1\rvert = \lvert x\generateur{x}\rvert, w_1$ is neither a factor of $\generateur{x}^{e_1}\prefix{x}$ nor of $\generateur{x}^{e_2}$.\\
Similarly for $w_2$, $n_{r'_s}(w_2) = n_{r'_s}(\generateur{x}) + 1$ and $n_{r_s}(\generateur{x}) = n_{r_s}(w_2) + 1$ and, by corollary \ref{number}, $w_2$ is not a conjugate of $\generateur{x}$, and because $\lvert w_2\rvert = \lvert x \rvert, w_2$ is neither a factor of $\generateur{x}^{e_1}\prefix{x}$ nor of $\generateur{x}^{e_2}$. \qed
\end{proof}

\begin{example}
Consider again $\generateur{x} = aaabaaaaaabaaaa$, $\prefix{x} = aaabaaaaaabaaa$ and $\suffix{x}=a$. We have $\lvert \generateur{x} \rvert = 15$, and:
\[ \generateur{x}\prefix{x}\generateur{x} = aaabaaaaaabaaaaaaa\rlap{$\overbrace{\phantom{baaaaaa\mathbf{baaaaaab}}}^{w_1}$}baaaaaa\underbrace{\mathbf{baaaaaab}aaaaaab}_{w_2}aaaa \]
The core of the interrupt is presented in bold.\\
The two factors $w_1$ and $w_2 = w_1 = baaaaaabaaaaaab$ (note that $w_2$ needs not be equal to $w_1$), starting at different positions, are not factors of $\generateur{x}^2$. Yet, the factor $aaaaaabaaaaaabaaaaaa$ of length $\lvert \generateur{x} \rvert + \lvert \lcs (\generateur{x}, \tilde{\generateur{x}}) \rvert + \lvert \lcp (\generateur{x}, \tilde{\generateur{x}}) \rvert$ and which contains the core of the interrupt is a factor of $\generateur{x}^2$. The same goes for the factors of length $\lvert \generateur{x} \rvert -1$ that starts and ends with the core of the interrupt, $aaaaaabaaaaaab$ and $baaaaaabaaaaaa$: they both are factors of $\generateur{x}^2$. For those reasons, the theorem can be regarded as tight
\end{example}

\section{Conclusion}

The key features of the core of the interrupt was understood while studying double squares. Ilie \cite{I05} provided an alternate and  shorter proof of Crochemore and Rytter's three squares lemma \cite{CR95}. We offer another concise proof within the framework of the core of the interrupt.

\begin{lemma}
In a word, no more that two squares can have their last occurrence starting at the same position.
\end{lemma}

\begin{proof}
Suppose that three squares $u_1^2, u_2^2, u_3^2, \lvert u_1 \rvert < \lvert u_2 \rvert < \lvert u_3 \rvert$ start at the same position. Because $u_2^2$ and $u_3^2$ start at the same position, we can write $u_2=\generateur{x}^{e_1}\prefix{x}, u_3=\generateur{x}^{e_1}\prefix{x}\generateur{x}^{e_2}$ for $\generateur{x} = \prefix{x}\suffix{x}$ a primitive word, $\prefix{x}$ a proper prefix of $\generateur{x}$ and $e_1 \geq e_2 \geq 1$, hence $u_3$ contains a core of the interrupt. Now, by synchronization principle, Theorem \ref{fw}, $u_1, \lvert u_1 \rvert < \lvert u_2 \rvert$, cannot end in the suffix $\lcs (\prefix{x}\suffix{x}, \suffix{x}\prefix{x})$ of $u_2$ (since $u_1$ has $\generateur{x}$ as a prefix) and ends before the core of the interrupt of $u_3$, but if $\lvert u_1^2 \rvert \geq \lvert u_3 \rvert$, the second occurrence of $u_1$ contains the core of the interrupt and a word of length $\lvert \generateur{x} \rvert$ that starts with it, while the first occurrence doesn't: which, by Theorem \ref{ath}, is a contradiction.
\end{proof}
\subsubsection{Thanks}

 to my supervisors Antoine Deza and Franya Franek for helpful discussions
and advices and to Alice Heliou for proof reading of a preliminary version of this article.

\end{document}